\documentclass[12pt]{article}
\usepackage{amssymb,latexsym,
amsfonts,amsbsy,curves,bezier,epic,color,eepic}
\usepackage[dvips]{graphics}

\newcommand{\qed}{\hbox{\rule[-2pt]
{3pt}{6pt}}}

\newtheorem{theo}{Theorem}[section]
\newtheorem{lem}[theo]{Lemma}
\newtheorem{pro}[theo]{Proposition}

\makeatletter

\@addtoreset{equation}{section}
\makeatother



\begin{title}
{Some Bounds for the Number of Blocks III}
\end{title}
\author{Etsuko Bannai and Ryuzaburo Noda}
\date{}
\begin{document}
\maketitle
\begin{abstract}
Let $\mathcal D=(\Omega, \mathcal B)$ be a pair of $v$ point set
$\Omega$ and a set $\mathcal B$ consists of 
$k$ point subsets of $\Omega$ which are called blocks.
Let $d$ be the
maximal cardinality of the intersections
between the
distinct two blocks in $\mathcal B$. The triple $(v,k,d)$  
is called the parameter of $\mathcal B$. 
Let $b$ be the number of the blocks in $\mathcal B$.
It is shown that inequality {\scriptsize
${v\choose d+2i-1}\geq b\left\{
{k\choose d+2i-1}
+{k\choose d+2i-2}{v-k\choose 1}+\cdots\right.$
$\left.+{k\choose d+i}{v-k\choose i-1}
\right\}$} holds for each $i$ satisfying $1\leq i\leq k-d$,
in the paper ``Some Bounds for the Number of Blocks",￼ Europ. J. Combinatorics 22 (2001), 91--94, by R. Noda. 

If $b$ achieves the upper bound for some 
$i$, $1\leq i\leq k-d$, then $\mathcal D$
is called a $\beta(i)$ design. 
In the paper mentioned above, an upper bound and a lower bound, {\scriptsize
$ \frac{(d+2i)(k-d)}{i}\leq  v \leq \frac{(d+2(i-1))(k-d)}{i-1}  $},
for $v$ of $\beta(i)$ design $\mathcal D$ are given.
In this paper we consider 
the cases when $v$ does not 
achieve the upper bound or lower bound given above, 
and get new more strict bounds for $v$ respectively.
We apply this bound to the problem of the perfect $e$-codes in the Johnson scheme, and improve the bound given by 
Roos in the paper ``A note on the existence of perfect constant 
weight codes'', Discrete Math. 47 (1983), 121--123.
\end{abstract}


 code


\section{Introduction}
This paper is a sequel to \cite{Noda-2001I} and \cite{Noda-2001II}. 
Let $\Omega$ be a $v$ point set. We denote by
${\Omega \choose k}$ the family of all the $k$ point subsets of $\Omega$.  
A {\it design} is a pair $\mathcal D=(\Omega,\mathcal B)$, where $\mathcal B$ is 
a subset of 
${\Omega\choose k}$.
The elements in $\mathcal B$ are called {\it blocks}.
We denote the number of the blocks in $\mathcal B$ by $b$.
A design $\mathcal D$ is called a $t$-$(v, k,\lambda)$ design 
if for every $t$ point subset of $\Omega$ there exist exactly $\lambda$ blocks containing it. 
A $t$-$(v, k,1)$ design is called a {\it Steiner system} and denoted by $S(t,k,v)$. 
A $t$-$(v,k,\lambda)$ design or a Steiner system S(t,k,v) is called 
{\it trivial} if $t=k$. In this paper we say that a design 
$\mathcal D=(\Omega,\mathcal B)$, or the block set $\mathcal B$, 
has the {\it parameter} $(v,k,d)$ 
if $\max\{|B\cap C| \mid B,C\in \mathcal B, B\neq C\}=d$. 
The number $d$ plays an essential role in our paper. 
We assume $d < k < v$ throughout. \\

The following is a slightly revised version of a basic Proposition given in
\cite{Noda-2001I} for designs with the parameter $(v,k,d)$.
\\

\noindent
{\bf Proposition 0}~{\bf (\cite{Noda-2001I}, see also
\cite{Cameron})}\\
{\it
Let $\mathcal D = (\Omega,\mathcal B)$ be a design 
and $(v,k,d)$ be the parameter of $\mathcal B$. Let $b=|\mathcal B|$.
Then the following inequality 
{\small
\begin{eqnarray}&&
{v\choose d+2i-1}\geq b\bigg\{
{k\choose d+2i-1}
+{k\choose d+2i-2}{v-k\choose 1}+\cdots+{k\choose d+i}{v-k\choose i-1}\bigg\}
\nonumber\\
&&\label{equ:1.1}
\end{eqnarray}}
holds for any $i$ satisfying $1\leq i\leq \frac{v-d+1}{2}$.
Moreover, equality in {\rm(\ref{equ:1.1})} holds for some $i$,
$1\leq i\leq \frac{v-d+1}{2}$, if and only if for 
any $(d+2i-1)$ point subset $X$ of $\Omega$, there exists a
block $B\in\mathcal B$ with $|X\cap B|\geq d+i$.\\
}

Proposition 0 follows by counting in two ways the cardinality of the set
$\{(X,B)\mid X\subset \Omega,~|X|=d+2i-1,~B\in\mathcal B,~|X\cap B|\geq d+i\}$. Note that for any $(d+2i-1)$ point subset $X$, there exists at most one block $B$ which satisfies $|X\cap B|\geq d+i$.\\

{\noindent
{\bf Definition 1}~{\bf ($\boldsymbol{\beta(i)}$ design)}\\
{\it We say that a design with the parameter
$(v,k,d)$ is a $\beta(i)$ design if the bound in
{\rm(\ref{equ:1.1})} is achieved for some $i$, with
$1\leq i\leq \frac{v-d+1}{2}$.
}\\

Obviously if the bound in (\ref{equ:1.1}) 
is achieved for some $i$ ,with $1\leq i\leq \frac{v-d+1}{2}$,
then we must have 
\begin{eqnarray} i\leq k-d.
\label{equ:1.2}
\end{eqnarray}  
A design $\mathcal D=(\Omega,\mathcal B)$
with the parameter $(v,k,d)$
is a $\beta(i)$ design if and only if for every $(d+2i-1)$ point 
subset $X\in {\Omega\choose d+2i-1}$ there exists a unique block $B\in \mathcal B$
with $|B\cap X|\geq d+i$.
In other words, a design $\mathcal  D$ is a $\beta(i)$ design 
if and only if ${\Omega\choose d+2i-1}$ has the
following partition parameterized by blocks in $\mathcal B$.
\begin{eqnarray}{\Omega\choose d+2i-1}=\bigcup_{B\in \mathcal B}
\left\{X\in {\Omega\choose d+2i-1}\mid
|X\cap B|\geq d+i\right\}.
\label{equ:1.3}
\end{eqnarray}
By Theorem 1 in \cite{Noda-2001I} the {\it complementary design} of a $\beta(i)$
design with the parameter $(v,k,d)$ is a $\beta(i')$ design with
\begin{eqnarray}
i'=k-d-i+1.
\label{equ:1.4}
\end{eqnarray}
The parameter $(v,k',d')$ of the complementary design of a $\beta(i)$
design with the parameter $(v,k,d)$ is given by $k'= v-k$ and $d'= v-2k+d$. 
Hence it follows that              
\begin{eqnarray}k'-d'=k-d.
\label{equ:1.5}
\end{eqnarray}
Steiner systems $S(t,k,v)$ are $\beta(1)$ designs (with $d=t-1$)
and conversely.
By (\ref{equ:1.4}) the complementary designs of Steiner systems $S(d+1,k,v)$
are $\beta(k-d)$ designs and conversely. We call a $\beta(i)$
design {\it trivial} if $i=k-d$ holds.
Only two non trivial $\beta(i)$ designs are known, the Steiner systems 
$S(5,8,24)$ $(i=2)$ and the complementary design of it ($i=3$). 
These are also a 
$\beta(1)$ design and a $\beta(4)$ design respectively.\\
\newpage
\noindent 
 {\bf Remark}{\it
 \begin{enumerate}
 \item If $\mathcal D$ is a $\beta(i)$ design with $i\geq 2$ which has 
 the parameter $(v,k,d=0)$, then $\mathcal D$ is a Steiner system
 $S(1,k,2k)$ with $k>1$. In this case $\mathcal D$ is a
$\beta(i)$ design for every $i$ with $ i\leq k$.
 \item If $\mathcal D$ is a $\beta(i)$ design with the parameter 
 $(v,k,d=k-1)$, then $i=1$ and $\mathcal D$ is a trivial Steiner system 
 $S(k,k,v)$.
 \end{enumerate}}
\noindent
 Proof of the Remark above is straightforward and omitted.\\

 In \cite{Noda-2001I}, Theorem 2 (1) the following inequalities for a 
 $\beta(i)$ design $\mathcal D$ are given.
 \begin{eqnarray} \frac{(d+2i)(k-d)}{i}\leq  v \leq \frac{(d+2(i-1))(k-d)}{i-1}.
 \label{equ:1.6}
 \end{eqnarray}
 Here $v$ attains the upper bound in (\ref{equ:1.6}) if and only if 
 $\mathcal D$ is a $\beta(i-1)$ design
and $v$ attains the lower bound in (\ref{equ:1.6}) if and only if $\mathcal D$ 
 is a $\beta(i+1)$ design.
Therefore if $d>0$, then a $\beta(i)$ design 
can not be a $\beta(j)$ design at the same time for $j$ with $|j-i|\geq 2$. 
The Steiner system $S(5,8,24)$ achieves the lower bound of 
(\ref{equ:1.6}) with $i=1$ and the complementary design of 
it achieves the upper bound of (\ref{equ:1.6}) with $i=4$. 
We remark that if a $\beta(i)$ design $\mathcal D$ achieves the upper bound 
of (\ref{equ:1.6}) then the complementary design of $\mathcal D$ 
achieves the lower bound of (\ref{equ:1.6}) and conversely.

In this paper we first give alternative bounds on the number of points in 
$\beta(i)$ 
designs if the upper or the lower bound given in (\ref{equ:1.6}) 
is not achieved (Theorems 1 and 2 in \S 2).
In \S 3 by making use of Theorems 1 and 2 we give new bounds on $v$ 
for a perfect $e$-code in the Johnson scheme $J(v,k)$ which improve the bound of Roos \cite{Roos} (Theorem 3). In  \S 4 we show that $k-d$ for a $\beta(i)$ design with parameter
  $(v,k,d)$ is bounded in terms of linear expressions of $i$ under some assumptions on $k-d$ and $i$ (Theorems 4 and 5).
                    
Our theorems are derived from the following basic propositions. \\

\noindent
{\bf Proposition 6}\\
{\it Let $\mathcal D$ be a $\beta (i)$ design with the parameter 
$(v,k,d)$, and let $c=k-d$. 
Assume that $i \geq 2$ and $v$ does not achieve the upper bound 
in {\rm(\ref{equ:1.6})}. 
Then we have
$g(v-k,d,c,i)\geq 0$ where
$g(x,d,c,i)$ is defined by
\begin{eqnarray}
&&g(x,d,c,i)=
\nonumber\\
&&(i-1)(i-2)x^2-(i-1)\bigg(2(c-i+1)d+2c(i-1)-3i+4\bigg)x
\nonumber\\
&&
+(c-i+1)(c-i+2)d^2+(c-i+1)\bigg((2i-3)c+3i-4\bigg)d
\nonumber\\
&&
+(i-1)(i-2)c^2+(i-1)(3i-4)c-2(i-1)^2(2i- 3).
\label{equ:1.7}
\end{eqnarray}
Moreover, $g(v-k,d,c,i)= 0$ if and only if $|\{B\in\mathcal B\mid |X\cap B|=d+i-2\}|$ does not 
depend on the choice of $X\in S_1$,
where $S_1$ is defined by 
\begin{eqnarray}
&& S_1=\{ X\subset
 \Omega \mid  |X|=d+2(i-2),~ |X\cap B|
\leq d+i-2, ~\forall B\in\mathcal B\}.
\label{equ:1.8}
\end{eqnarray}
}

\noindent
{\bf Proposition 7}\\
{\it
Let $\mathcal D$ be a $\beta (i)$ design with the parameter 
$(v,k,d)$. Assume that $k-d\geq i+1$ and $v$ does not achieve the lower
bound in {\rm(\ref{equ:1.6})}. Then 
$h(v-k,d,c,i)\geq 0$ where $h(x,d,c,i)$
is defined by
\begin{eqnarray}
&&h(x,d,c,i)=i(i+1)x^2-i\bigg(2(c-i)d+2ic+3i+1\bigg)x
\nonumber\\
&&
+(c-i)(c-i-1)d^2+(c-i)\bigg((2i+1)c-3i-1\bigg)d
\nonumber\\
&&
+i(i+1)c^2-i(3i+1)c+2i^2(2i+1).
\label{equ:1.9}
\end{eqnarray}
Moreover, $h(v-k,d,c,i))=0$
if and only if $|\{B\in \mathcal B \mid |X\cap B|=d+i+1\}|$
does not depend on the choice of $X\in S_2$,
where
$S_2$ is defined by 
\begin{eqnarray}
&& S_2 = \{ X \subset \Omega \mid |X|=d+2(i+1),~|X\cap B|
\leq d+i+1,~\forall B\in\mathcal B\}.
\label{equ:1.10}
\end{eqnarray}
}

\noindent
It is shown in Lemma \ref{lem:2.1} that $S_1$ and $S_2$ are not empty under the assumptions of Propositions 6 and 7 respectively.
\noindent
Proofs of Propositions 6 and 7 are very long and
involved, so we give them later in \S 5.
In \S 6 we present some open problems 
in $\beta(i)$ designs.\\

Finally we would like to announce the following.\\

\noindent
{\bf Theorem}\quad {\it Let $\mathcal D$ be a $\beta(i)$ design with the 
parameter $(v,k,d)$. 
Then the following polynomial in $t$ of degree $2(i-1)$ has at least one positive integral zero and $\mathcal D$ is a $t$-$(v,k,\lambda)$ design for $t$ 
which is the smallest of positive integral zeros of it.
\begin{eqnarray}
\sum_{j=0}^{i-1}(-1)^j{t+1\choose j}
\sum_{s=0}^{i-j-1}{k-t+j-1\choose d-t+2i-2-s}{v-k-j\choose s}.
\label{equ:1.11}
\end{eqnarray}
}

The proof of the above Theorem is given in a forthcoming paper \cite{Noda}.\\

\noindent
{\bf Remark}\\
{\it The inequality given in Proposition 0 is also known by the 
researchers who studied
maximal intersecting systems of finite sets which is originated by
Erd\"{o}s, Ko and Rado 
in {\rm1961~\cite{Erdos-C-K-1961}}. 
Ahlswede-Aydinian-Khachatrian in {\rm 2001 \cite{Ahlswede-A-K-2001}}
proved that the existence of the partition of 
 ${\Omega\choose d+2i-1}$ given above is equivalent 
 to an existence of the $D$-diameter perfect code with diameter $D=k-d+1$.
 So the existence of a $\beta(i)$ design with the parameter
 $(v,k,d)$ is equivalent to the existence of the $D$-diameter perfect code with diameter $D=k-d+1$. 
They obtained upper and lower bound of the cardinality of a
$D$-diameter perfect code which are equivalent to the bounds given in
{\rm(\ref{equ:1.6})}.
We found out this fact very recently. 
For more information on intersecting systems of finite sets and 
perfect $D$-diameter codes please refer to
{\rm\cite{Ahlswede-K-1997}, \cite{Ahlswede-A-K-2001}} and {\rm\cite{Delsarte-1973}}. }

\section{Some bounds on $\boldsymbol{v=|\Omega|}$}
In this section we give alternative bounds on $v$ for $\beta(i)$ designs 
with the parameter $(v,k,d)$
which do not achieve the bounds in (\ref{equ:1.6}). 
For this purpose we use the following two 
families of point subsets, $S_1$ and $S_2$, defined in Proposition 6 (\ref{equ:1.8}) and Proposition 7 (\ref{equ:1.10}) 
respectively, 
for a $\beta(i)$ design $\mathcal D=(\Omega, \mathcal B)$ 
with the parameter $(v,k,d)$.
\begin{eqnarray}
&& S_1=\{ X\subset
 \Omega \mid  |X|=d+2(i-2),~ |X\cap B|
\leq d+i-2, ~\forall B\in\mathcal B\}.
\nonumber\\
&& S_2 = \{ X \subset \Omega \mid |X|=d+2(i+1),~|X\cap B|
\leq d+i+1,~\forall B\in\mathcal B\}.
\nonumber
\end{eqnarray}
First we prove the following.
\begin{lem}
 Let $\mathcal D$ be a $\beta(i)$ design having the parameter
$(v,k,d)$. Then 
 \begin{enumerate}
 \item $S_1$ is not empty if $i\geq 2$, and
\item	 $S_2$ is not empty if $k-d\geq i+1$.
 \end{enumerate}
 \label{lem:2.1}
\end{lem}
{\bf Proof}\\
Let $k',~d'$ and $i'$ denote the parameters of the 
complementary design $\mathcal D'$ of $\mathcal D$.
Then 
$$k' = v-k,\quad d'=v-2k+d \quad\mbox{and} \quad i'=k-d-i+1$$
by (\ref{equ:1.4}).

(1) If there exists a $d+2(i-2)$ point subset $X$ such that $|X\cap B|=d+i-2$ 
and $|X\cap(\Omega\backslash B)|=i-2$ for some block $B$, 
then $X$ belongs to $S_1$. 
Therefore in order to prove (1) it suffices to show that
$k\geq d+i-2$
and
$v-k\geq i-2\geq 0$.
We have $k-d\geq i-2$ 
by (\ref{equ:1.2}) and $i\geq 2$ by our assumption. 
Also since 
$$(v-k)-(i-2)=k'-(k'- d'-i'-1)=d' +i'+1>0,$$ 
we have $v-k > i-2$.\\
(2) Similarly if there exists a $d+2(i+1)$ point subset $Y$ 
such that $|Y\cap B|=d+i+1$ and $|Y\cap (\Omega\backslash B)|=i+1$ for some block $B$, then $Y$ belongs to $S_2$. We have $k\geq d+i+1$ by our assumption. Therefore in order to prove (2) it suffices to show that
$v-k\geq i+1$.
Since $d' \geq 0$ we have $v-k\geq  k-d\geq  i+1$.
\hfill\qed\\

We can now state our results.\\

 \noindent
{\bf Theorem 1}\quad
{\it Let $\mathcal D$ be a $\beta (i)$ design satisfying $i\geq 3$ and 
$(v,k,d)$ be the parameter of $\mathcal B$ with $v\geq 2k$.
Let $c=k-d$. Assume that $v$ does not achieve the upper bound in 
{\rm(\ref{equ:1.6})}.
 Then we have
\begin{eqnarray}
v\leq k+\gamma_1= c+d+\gamma_1,
\label{equ:2.1}
\end{eqnarray}
where $\gamma_1$ is given
by the following formula.
\begin{eqnarray}
&&\gamma_1=\frac{(2c-2i+2)d+2(i-1)c-3i+4}{2(i-2)}
-\frac{1}{2(i-1)(i-2)}\times
\nonumber\\
&&
\bigg\{4(i-1)(c-i+1)(c-2i+3)d^2+4(i-1)(3i-4)(c-i+1)(c-2i+3)d
\nonumber\\
&&
+(i-1)^2\bigg(4(2i-3)c^2-4(3i-4)(2i-3)c
+16i^3-63i^2+80i-32\bigg)\bigg\}^{\frac{1}{2}}.
\label{equ:2.2}
\end{eqnarray}
Moreover, equality in {\rm(\ref{equ:2.1})} holds
if and only if the number of the blocks $B$, for which $|B\cap X|
=d+i-2$, does not depend on the choice of 
$X\in S_1$.}
\\

\noindent
{\bf Theorem 2}\\
{\it Let $\mathcal D$ be a $\beta (i)$ design and $(v,k,d)$
be the
parameter of $\mathcal B$ satisfying 
$v \leq  2k$ and $k-d \geq  i+1$. 
Let $c=k-d$.
Assume that $v$ does not achieve the lower bound in {\rm(\ref{equ:1.6})}.
Then we have
\begin{eqnarray}
v\geq k+\gamma_2=c+d+\gamma_2,
\label{equ:2.3}
\end{eqnarray}
where $\gamma_2$ is given by the following formula.
\begin{eqnarray}
&&\gamma_2=\frac{2(c-i)d+2ic+3i+1}{2(i+1)}
+\frac{1}{2i(i+1)}\times
\nonumber\\
&&
\bigg\{
-4i(c-i)(c-2i-1)d^2-4i(3i+1)(c-i)(c-2i-1)d
\nonumber\\
&&
-i^2\bigg(4(2i+1)c^2-4(2i+1)(3i+1)c+16i^3+15i^2+2i-1\bigg)\bigg\}
^{\frac{1}{2}}.
\label{equ:2.4}
\end{eqnarray}
Moreover, equality in {\rm(\ref{equ:2.3})} holds 
if and only if the number of blocks
$B$, for which $|B\cap X|=d+i+1$, does not depend on the
choice of $X\in S_2$.}\\

In order to prove Theorems 1 and 2 we need the following.

\begin{lem}
 Let $\mathcal D$ be a $\beta(i)$ design having the parameter
$(v,k,d)$, and let $c=k-d$. Then we have the following.
 \begin{enumerate}
 \item If $v\geq 2k$, then $c\geq 2(i-1)$ with equality if and only if 
 $\mathcal D$ is a $\beta(i-1)$ design with $v=2k$.
\item	If $v\leq  2k$, then $c\leq  2i$ with equality if and only if $\mathcal D$
is a $\beta(i+1)$ design with $v=2k$.
\item	If $v=2k$, then $2(i-1)\leq c\leq 2i$ holds. 
 \end{enumerate}
 \label{lem:2.2}
\end{lem}
{\bf Proof}\\
(1) Assume that $v\geq  2k$. Then by (\ref{equ:1.6}),
  $$2k\leq v \leq  \frac{(d+2(i-1))c}{i-1}=\frac{dc}{i-1} +2c,$$
  hence
$2d \leq\frac{dc}{i-1}$. Therefore $2(i-1)\leq c$ with equality if and only 
  if $v=2k$ and the upper bound of (\ref{equ:1.6}) is attained.
 Recall that $\mathcal D$ is a $\beta(i-1)$ design 
  if and only if the upper bound of (\ref{equ:1.6})
is achieved.
  \\
(2) We can prove (2) in the same way as (1) by using the lower 
bound of (\ref{equ:1.6}).\\
(3)	An immediate consequence of (1) and (2).
\hfill\qed\\

Let 
\begin{eqnarray}&&x_1=\frac{d(c-i)}{i}+c,
\label{equ:2.5}\end{eqnarray}
and let
\begin{eqnarray}&& x_2=\frac{d(c-i+1)}{i-1}+c.
\label{equ:2.6}
\end{eqnarray}
 Then (\ref{equ:1.6}) implies
\begin{eqnarray}
x_1\leq v-k\leq x_2.
\label{equ:2.7}
\end{eqnarray}
For the proof of Theorem 1, we use the
upper bound on $v-k$, $x_2$
and for the proof of Theorem 2, we use the
lower bound on $v-k$, $x_1$ as
given above.\\

\noindent{\bf Proof of Theorem 1.}\\
By assumption, $\mathcal D$ is a $\beta(i)$ design and
is not a $\beta(i-1)$ design with the parameter set $(v,k,d)$
satisfying $i\geq 3$ and $v\geq 2k$. Then Proposition 6 implies that
$g(v-k,d,c,i)\geq 0$. We show that $g(x_2,d,c,i)<0$.
We have
\begin{eqnarray}
&&g(x_2,d,c,i)=(i-1)(i-2)\left(\frac{d(c-i+1)}{i-1}+c\right)^2
\nonumber\\
&&
-(i-1)\bigg(2(c-i+1)d+2c(i-1)-3i+4\bigg)\left(\frac{d(c-i+1)}{i-1}+c\right)
\nonumber\\
&&
+(c-i+1)(c-i+2)d^2+(c-i+1)\bigg((2i-3)c+3i-4\bigg)d
\nonumber\\
&&
+(i-1)(i-2)c^2+(i-1)(3i-4)c-2(i-1)^2(2i- 3)
\nonumber\\
&&
=-\frac{(d+i-1)(c-i+1)}{i-1} \bigg((c-2i+2)d+2(i-1)c-4i^2+10i-6\bigg).
\label{equ:2.8}
\end{eqnarray}
Since $v\geq 2k$ by our assumption, Lemma \ref{lem:2.2} (1)
implies that $c\geq 2(i-1)$. Therefore
\begin{eqnarray}
&&(c-2i+2)d+2(i-1)c-4i^2+10i-6
\nonumber\\
&&
\geq 2(i-1)(2(i-1))-4i^2+10i-6=2(i-1)>0.
\label{equ:2.9}
\end{eqnarray}
By (\ref{equ:2.8}) and (\ref{equ:2.9})
we have 
\begin{eqnarray}&&g(x_2,d,c,i)< 0.
\label{equ:2.10}
\end{eqnarray} 
Since $g(x,d,c,i)$ is
a polynomial in $x$ of degree $2$ 
and the coefficient of $x^2$
is positive, $g(x,d,c,i)$ has two distinct real zeros.
Let $\gamma_1$ be the smaller one.
Since $g(v-k,d,c,i)\geq 0$ and $v-k\leq x_2$, we must have
$v-k\leq \gamma_1$. 
Moreover the equality 
$v-k= \gamma_1$ holds if and only if $g(v-k,d,c,i)=0$. 
Hence, Proposition 6 completes the proof of 
Theorem 1.
\hfill\qed\\

\noindent{\bf Proof of Theorem 2.}\\
By assumption, $\mathcal D$ is a $\beta(i)$ design 
and is not a $\beta(i+1)$ design with
the parameter set $(v,k,d)$ satisfying $v\leq 2k$ and
$c=k-d\geq i+1$. Then Proposition 7 implies 
that $h(v-k,d,c,i)\geq 0$. We show that $h(x_1,d,c,i)<0$.
We have
\begin{eqnarray}
&&h(x_1,d,c,i)=i(i+1)\left(\frac{d(c-i)}{i}+c\right)^2-i\bigg(2(c-i)d+2ic+3i+1\bigg)\left(\frac{d(c-i)}{i}+c\right)
\nonumber\\
&&
+(c-i)(c-i-1)d^2+(c-i)\bigg((2i+1)c-3i-1\bigg)d
\nonumber\\
&&
+i(i+1)c^2-i(3i+1)c+2i^2(2i+1)
\nonumber\\
&&
=\frac{(d+i)(c-i)}{i}\bigg((d+2i)(c-2i)-2i\bigg).
\label{equ:2.11}
\end{eqnarray}
Since $v\leq 2k$ by our assumption, Lemma \ref{lem:2.2} (2) implies that 
$c\leq 2i$. Hence, it follows by (\ref{equ:2.11}) that 
\begin{eqnarray}
h(x_1,d,c,i)<0.
\label{equ:2.12}
\end{eqnarray}
Since $h(x,d,c,i)$ is a polynomial in $x$ of degree 2 and
the coefficient of $x^2$ is positive, $h(x,d,c,i)$ has two 
distinct real zeros. Let $\gamma_2$ be the larger one.
Since $h(v-k,d,c,i)\geq 0$ and $v-k\geq x_1$,
we must have $v-k\geq \gamma_2$.
Moreover the equality, $v-k= \gamma_2$, holds 
if and only if $h(v-k,d,c,i)= 0$.
Hence, Proposition 7
completes the proof of Theorem 2.
\hfill\qed\\

\noindent
{\bf Remark}{\it
\begin{enumerate}
\item If $\mathcal D$ achieves the upper bound in {\rm(\ref{equ:1.6})}, then for every $X\in S_1$ there exist
$\frac{v-d-2(i-2)}{k-d-i+2}$ blocks in $\mathcal B$ which have $d+i-2$ points in common with $X$.
\item If $\mathcal D$ achieves the lower bound  in {\rm(\ref{equ:1.6})}, then for every $X\in S_2$ there exist
$\frac{d+2(i+1)}{i+1}$ blocks in $\mathcal B$ which have $d+i+1$ points in common with $X$.
\end{enumerate}}

\section{Application to perfect $e$-codes}

The {\it Johnson graph} (or the {\it Johnson scheme}) $J(v,k)$ is a graph whose set of
vertices is ${\Omega\choose k}$ for a $v$ point set $\Omega$. 
Two vertices $x$ and $y$ are {\it adjacent} if and only if $|x\cap y|=k-1$. 
The {\it distance} $d(x,y)$ between two vertices $x$ and $y$ is the {\it length} of the {\it shortest path} which connects these vertices, that is , $d(x,y)=k-|x\cap y|$. 
A {\it code} $\mathcal C$ in $J(v,k)$ is a subset of ${\Omega\choose k}$.
The {\it minimum distance} of $\mathcal C$ is defined by 
$d(C)=\min\{d(x,y) |x\neq y,~  x,y\in \mathcal C\}$.
A code $\mathcal C$ is called a {\it perfect $e$-code} in $J(v,k)$ 
if the $e$-spheres with centers at the code words of $\mathcal C$ form a 
partition of ${\Omega\choose k}$. In other words, $\mathcal C$ is a 
perfect $e$-code if for each element $x\in{\Omega \choose k}$ 
there exists a unique element $c\in \mathcal C$ such that $d(x,c)\leq e$. 
Clearly, the minimum distance of a perfect $e$-code is $2e+1$.

Then it is easily seen that a $\beta(i)$ set 
$\mathcal B$ having the parameter $(v,k,d)$ with $k-d=2i-1$ 
is a perfect $(i-1)$-code in the Johnson scheme $J(v,k)$ 
and conversely.
 As for a perfect $e$-code in $J(v,k)$, Roos \cite{Roos} proved the inequality 
 $v\leq (k-1)(2e +1)/e$. This corresponds to the upper bound of (\ref{equ:1.6}). 
 T. Etzion and M. Schwartz (\cite{Etzion-S-2004}, Theorem 13) showed that the bound of Roos is not achievable. For the details of perfect $e$-codes in the 
Johnson scheme $J(v,k)$ please refer to \cite{Hauck-1982}.

By making use of the same argument as employed in the proofs of
Theorems 1 and 2 we prove the following theorem which improves 
the upper bound of Roos.\\

\noindent
{\bf Theorem 3}\\
{\it Let $e$ be an integer satisfying $e\geq 2$ and
$\mathcal C$ be a perfect $e$-code in the Johnson scheme $J(v,k)$. 
Then we have the following inequality.
{\small\begin{eqnarray}
&&\frac{ 2(e+1)k}{e+2}+\frac{7e+6}{2(e+2)}
+\frac{\sqrt{A_2}}{2(e+1)(e+2)}\leq v
\leq\frac{2ke}{e-1} -\frac{7e+1}{2(e-1)} - \frac{\sqrt{A_1}}{2e(e-1)},
\nonumber\\
&&  \label{equ:3.1}
\\
&\mbox{where}&\nonumber\\
&&A_1=e\left\{8(e+1)\left(k-\frac{e+3}2\right)^2-(e+2)(e-1)^2\right\}\nonumber\\
&\mbox{and}&\nonumber\\
 &&A_2=(e+1)\left\{8e\left(k-\frac{e-2}{2}\right)^2-(e-1)(e+2)^2)\right\}.
\nonumber\end{eqnarray}
}\\
Moreover, $v$ achieves the upper bound in {\rm(\ref{equ:3.1})} 
if and only if for any word $y$ of length $v$ 
and of weight $k-3$, for which every
$u\in\mathcal C$ satisfies
 $|y\cap u|\leq k-e-2$,
the number of $u\in\mathcal C$ 
which satisfies $|u\cap y|= k-e-2$, is invariant,
that is, independent of the choice of such $y$.
Also, the lower bound holds if and only if for every word 
$y$ of length $v$ and of weight $k+3$, for which 
every $u\in \mathcal C$ satisfies $|y\cap u|\leq k-e+1$,
the number of $u\in\mathcal C$, which satisfies $|u\cap x|=k-e+1$, 
is invariant.}\\

\noindent
{\bf Proof}\\
Perfect $e$-codes in the Johnson scheme $J(v,k)$ 
can be regarded as block set of
$\beta(i)$ designs having the parameter
$(v,k,d)$, with $i=e+1$ and $d=k-2e-1=k-2i+1$.
Then since $2(i-1)< c=k-d< 2i$, we have
$g(x_2,d,c,i)<0$ by (\ref{equ:2.10}) and $h(x_1,d,c,i)<0$ by (\ref{equ:2.12}).
Hence the same argument as given in the proofs
of Theorems 1 and 2 yields Theorem 3.
\hfill\qed\\

\section{Bounds on $k-d$}
In this section we give bounds on $k-d$
of a $\beta(i)$ design $\mathcal D$ in terms of $i$.
These bounds are very
important and useful for the classification of $\beta(i)$ designs.\\

\noindent
{\bf Theorem 4}\\
{\it Let $\mathcal D$ be a $\beta(i)$ design with $i \geq 3$ 
and let $(v,k,d)$ be the parameter of $\mathcal D$,
where
$v\geq 2k$. Then the following hold.
\begin{enumerate}
\item If $\mathcal D$ is not a $\beta(i-1)$ design,
then 
\begin{eqnarray}2(i-1)\leq k-d < \frac{i( 3i-2 + \sqrt{i^2+12i-12})}{2(i-2)}.
\end{eqnarray}
In particular
\begin{eqnarray}2(i-1)\leq  k-d\leq  2i + 6
\end{eqnarray} for $i\geq  8$.
\item If $\mathcal D$ is a $\beta(i-1)$ design with $d>0$, then
\begin{eqnarray}
2(i-1)\leq k-d \leq 2i+4
\end{eqnarray}
for $i \geq 9$.
\end{enumerate}}

\noindent
{\bf Proof}\\
(1) Let $g(x,d,c,i)$ be the polynomial defined in 
Proposition 6. Let $x_1=\frac{d(c-i)}{i}+c$ and 
$x_2=\frac{d(c-i+1)}{i-1}+c$ be the lower bound
and the upper bound  
on $v-k$
defined in (\ref{equ:2.5}) and (\ref{equ:2.6}) respectively. 
We have $g(x_2,d,c,i)<0$ by (2.10) and $g(v-k,d,c,i)\geq 0$
by Proposition 6. Then since 
the coefficient of $x^2$
in the 
quadratic polynomial $g(x,v,k,d)$ is positive and
$x_1\leq v-k$, it follows that  $g(x_1,d,c,i)\geq 0.$
By substituting $x_1$ in $g(x, d, c, i)$ we have
\begin{eqnarray}
&&\frac{1}{i^2}~g(x_1,d,c,i)
=-\bigg((i-2)c^2-(3i^2-2i)c+2i^2(i-1)\bigg)d^2
\nonumber\\
&&
-i\bigg((3i-2)c-3i^2+4i\bigg)\bigg(c-2(i-1)\bigg)d-2i^2(i-1)(c-i+1)(c-2i+3),
\nonumber\\
&&c-2i+3\geq 0.
\label{equ:N4.4}
\end{eqnarray}
Since $i>2$ and $c\geq 2(i-1)$, 
it follows that
\begin{eqnarray}
&&(3i-2)c-3i^2+4i\geq (3i-2)(2i-2)-3i^2+4i=3i(i-2)+4>0,\nonumber\\
&&c-2(i-1)\geq 0,\nonumber\\
&&2i^2(i-1)(c-i+1)(c-2i+3)\geq 2i^2(i-1)^2>0.
\label{equ:N4.5}
\end{eqnarray}
By (\ref{equ:N4.4}) and (\ref{equ:N4.5}), 
we have
\begin{eqnarray}
&&(i-2)c^2-(3i^2-2i)c+2i^2(i-1)< 0.
\label{equ:4.6}
\end{eqnarray}
This implies 
\begin{eqnarray}
&&2(i-1)\leq c<\frac{i(3i-2+\sqrt{i^2+12i-12})}{2(i-2)}
\end{eqnarray}
It is easy to see that if $i\geq 8$, then
$\frac{i(3i-2+\sqrt{i^2+12i-12})}{2(i-2)}<2i+7$ holds.
This completes the proof of (1).\\
(2) As mentioned in the Introduction, if $\mathcal D$ is a $\beta(i)$ design
and a $\beta(i-1)$ design 
with $d>0$, then
$\mathcal D$ cannot be a $\beta(i-2)$ design. 
Therefore we can apply (1) to (2) and we 
have $c\leq 2(i-1)+6$ for $i-1\geq 8$.
This completes the proof of Theorem 4.
\hfill\qed\\

\noindent
{\bf Theorem 5}\\
{\it Let $\mathcal D$ be a $\beta(i)$ design with the 
parameter $(v,k,d)$.
Assume that 
$v\leq 2k$ and $c=k-d\geq i+2$. Then the following hold.
\begin{enumerate}
\item If $\mathcal D$ is not a $\beta(i+1)$ design and if $i\geq 18$,
then
\begin{eqnarray}
\frac{(i-1)(3i-1+\sqrt{i^2-14i+1})}{2(i+1)} < k-d\leq 2i.
\label{equ:4.8}
\end{eqnarray}
In particular 
$$2i-8 \leq  k-d \leq 2i.$$
\item If $\mathcal D$ is not a $\beta(i+1)$ design
and if $c=k-d\geq 3$, then
\begin{eqnarray}
\frac{(i-1)(3i-1+\sqrt{i^2-14i+1})}{2(i+1)} < k-d\leq 2i,
\end{eqnarray}
for $i\geq 15$. In particular 
$$2i-8\leq k-d\leq 2i,$$
if $i\geq 16$.
\item If $\mathcal D$ is a $\beta(i+1)$ design with $d>0$ and if $i\geq 17$, 
then 
\begin{eqnarray}2i-6\leq k-d \leq 2i.
\end{eqnarray}
\end{enumerate}}

\noindent
{\bf Proof}\\
Our assumption $v\leq 2k$ and Lemma \ref{lem:2.2} imply that $c\leq 2i$.
Let $\mathcal D'$ be the complementary design of $\mathcal D$.
Then $\mathcal D'$ is a $\beta(i')$ design with
the parameter $(v,k',d')$ where 
$$i'=k-d-i+1=c-i+1,\quad k'=v-k, \quad d'=v-2k+d$$
 (see (\ref{equ:1.4})).  
Then 
$v\leq 2k$ and $k-d\geq i+2 $ imply 
$$v\geq 2k'\quad\mbox{and}\quad i'=k-d-i+1\geq 3.$$ 
Therefore we can apply 
Theorem 2 to $\mathcal D'$.\\
(1) If $\mathcal D$ is not a $\beta(i+1)$ design, then 
$\mathcal D'$ is not a $\beta(k-d-(i+1)+1)$ design,
i.e. $\mathcal D'$ is a $\beta(i')$, but not a
 $\beta(i'-1)$ design. Hence we can use the inequality 
 (\ref{equ:4.6}) for $\mathcal D'$.
 By substituting $i'=c-i+1$ and $c'=c$,
we obtain
$$(i'-2)c^2-3(i'^2-2i')c+2i'^2(i'-1)<0$$ 
and so
$$(c-i-1) c^2-3((c-i+1)^2-2(c-i+1))c+2(c-i+1)^2(c-i)<0,$$
hence
\begin{eqnarray}
(i+1)c^2-(3i-1)(i-1)c+2i(i-1)^2>0.
\label{equ:4.11}
\end{eqnarray}
By (\ref{equ:4.11}), for $i\geq 14$ we have
$$c>\frac{(i-1)(3i-1+\sqrt{i^2-14i+1})}{2(i+1)}\quad \mbox{or}\quad
c<\frac{(i-1)(3i-1-\sqrt{i^2-14i+1})}{2(i+1)}.$$
Assume that the latter holds. 
Then by our assumption on $c$ we have
$$ i+2\leq c <\frac{(i-1)(3i-1-\sqrt{i^2-14i+1})}{2(i+1)},$$ 
and so
$$(i-1)(\sqrt{i^2-14i+1}) <i^2 -10i -3,$$
hence
$$(i+1)(i^2-17i-2)<0.$$
This implies $i\leq 17$. Therefore if $i\geq 18$, 
we have
$c>\frac{(i-1)(3i-1+\sqrt{i^2-14i+1})}{2(i+1)}$,
hence $c \geq 2i-8$.
\\
(2) As in the above the inequality 
$i+3 \leq \frac{(i-1)(3i-1-\sqrt{i^2-14i+1})}{2(i+1)}$
gives $(i+1)(i^2-14i-3)<0$. This implies $i\leq 14$. 
Therefore if  $i\geq 15$ we have
$c>\frac{(i-1)(3i-1+\sqrt{i^2-14i+1})}{2(i+1)}$ and 
hence we have $c \geq 2i-8$ if $i\geq 16$.\\
(3) If D is a $\beta(i+1)$ design then $\mathcal D$ 
is not a $\beta(i+2)$ design. Hence by (1)
we have $2(i+1)-8\leq c\leq 2i$ for $i+1\geq 18$.
 This completes the proof of Theorem 5.
 \hfill\qed\\

\section{Proof of Basic  Propositions}
In this section we give the proofs of the basic Propositions 6 and 7 
which are stated in
\S 1 and already used in the proofs of our main 
results in the previous sections.

The proof needs lots of preparations.
In the statements of Propositions 6 and 7,
we defined two families of subsets of $\Omega$, $S_1$ and $S_2$.
We state the definitions of $S_1$ and $S_2$ here again.
\begin{eqnarray}
&& S_1=\{ X\subset
 \Omega \mid  |X|=d+2(i-2),~ |X\cap B|
\leq d+i-2, ~\forall B\in\mathcal B\}.
\nonumber\\
&& S_2 = \{ X \subset \Omega \mid |X|=d+2(i+1),~|X\cap B|
\leq d+i+1,~\forall B\in\mathcal B\}.\nonumber
\end{eqnarray}
By Lemma \ref{lem:2.1}, $S_1$ and $S_2$ are not empty 
under the assumption of Propositions 6 and 7.

We first introduce the following 
combinatorial formula which we use later in the proof of
Proposition \ref{pro:5.5}.

\begin{lem}
Let $m,~s$ and $j$ be integers satisfying $m>0$ and
$s\geq j>0$. Then the following holds.
\begin{eqnarray}
{m\choose s}+{m\choose s-1}
{j\choose 1}+{m\choose s-2}{j\choose 2}+\cdots+{m\choose s-j}{j\choose j}= \frac{{m+j\choose j}{m\choose s-j}}{{s\choose j}}.
\label{equ:5.1}
\end{eqnarray}
\label{lem:5.1}
\end{lem}
We will give the proof of Lemma \ref{lem:5.1} later, at the end of
this section.\\


Now we are ready to start the proofs of Propositions 6 
and 7.
We first prove the following.
\begin{pro}
Let $\mathcal D=(\Omega,\mathcal B)$ be a $\beta(i)$ design
with the parameter $(v,k,d)$ and $B\in \mathcal B$ arbitrarily fixed. Let
$\mu_d=|\{C\in\mathcal B\mid |C\cap B|=d\}$ and $c=k-d$.
Then 
\begin{equation}\mu_d=\frac{{v-k\choose i}{k\choose d+i-1}}{{c\choose i}{c\choose i-1}}
\label{equ:5.2}
\end{equation}
holds. In particular $\mu_d$ is independent of the choice of $B\in\mathcal B$. 
\label{pro:5.2}
\end{pro}
{\bf Proof}\\
We count the cardinality of the following set in two different ways.
\begin{eqnarray}
\left\{ (X,C)~\bigg| \begin{array}{l}X \subset \Omega, C\in\mathcal B,
|X\cap B|=d+i-1,~|X\cap(\Omega\backslash B)|=i, \\
|X\cap B\cap C|=d,
   |X\cap(\Omega\backslash B) \cap C|=i
   \end{array}\right\}.
   \label{equ:5.3}
\end{eqnarray}
Let $X\subset \Omega$ be any fixed subset of $\Omega$ satisfying
 $|X\cap B|=d+i-1$ and $|X\cap (\Omega\backslash B)|=i$.
Since $\mathcal D$ is a $\beta(i)$ design, there exists 
a unique block $C\in \mathcal B$ satisfying 
$|X\cap C|\geq  d+i$
(see Definition 1 and (\ref{equ:1.3})).
Then since $|B\cap C|\leq d$, we must 
have $|X\cap B\cap C|= d$ and $|X\cap(\Omega\backslash B)\cap C| =i$. 
Hence the cardinality of the set defined in (\ref{equ:5.3}) equals
${k\choose d+i-1}{v-k\choose i}$.
On the other hand, for any block $C$ satisfying
$|B\cap C|= d$, the number of subsets $X\subset\Omega$ satisfying
$|X\cap B|=d+i-1$,  $|X\cap B\cap C|=d$ and 
$|X\cap(\Omega\backslash B)\cap C|= i$ equals
${k-d\choose i-1}{k-d\choose i}$. 
Hence the cardinality of the set defined in (\ref{equ:5.3})
equals $\mu_d{k-d\choose i-1}{k-d\choose i}$.
This completes the proof of Proposition
\ref{pro:5.2}.
\hfill\qed\\
\begin{pro}
Let $\mathcal D=(\Omega,\mathcal B)$ be a $\beta(i)$ design
with the parameter $(v,k,d)$.
Let $i\geq 2$, $c=k-d$, $b=|\mathcal B|$ and $n=|S_1|$, where $S_1$ is as defined in 
(\ref{equ:1.8}).
Then we have
\begin{eqnarray}
\frac{b}{n}{c+d \choose d+i-2}{v-k\choose i-2}
\leq 1+\frac{(i-1)(v-k-i+1)(v-k-i+2)}{(d+i-1)(c-i+1)(c-i+2)}.
\label{equ:5.4}
\end{eqnarray}
Moreover, equality in {\rm(\ref{equ:5.4})} holds
if and only if $|\{B\mid |B\cap X|=d+i-2\}|$ does not depend on the choice of $X\in S_1$. 
\label{pro:5.3}
\end{pro}
{\bf Proof}\\
For $X\in S_1$ let $\alpha_X$ be defined as follows.
\begin{eqnarray}
\alpha_X=|\{B\in \mathcal B\mid |X\cap B|=d+i-2 \}|.
\nonumber
\end{eqnarray}
Let $P=\sum_{X\in S_1}\alpha_X$ and 
$Q=\sum_{X\in S_1}\alpha_X(\alpha_X-1)$.
Then
{\small\begin{eqnarray}
&&0\leq \sum_{X\in S_1}\left(\alpha_X-\frac{P}{n}\right)^2
=\sum_{X\in S_1}\bigg(\alpha_X^2-\frac{2P}{n}\alpha_X+\frac{P^2}{n^2}\bigg)
=\sum_{X\in S_1}\alpha_X^2-\frac{P^2}{n}=Q+P-\frac{P^2}{n}
\nonumber
\end{eqnarray}}
implies 
\begin{eqnarray}
\frac{P}{n}\leq 1+\frac{Q}{P},
\label{equ:5.5N}
\end{eqnarray}
where equality holds if and only if $\alpha_X=\frac{P}{n}$
(constant),
i.e. $|\{B\mid |B\cap X|=d+i-2\}|$ does not depend on the choice of $X\in S_1$.
We remark that if a $\beta(i)$ design $\mathcal D$ is 
also a $\beta(i-1)$ design, 
then $\alpha_X$ is a constant given by
$\alpha_X=\frac{v-d-2(i-2)}{k-d-i+2}=\frac{v-k+c-2(i-2)}{c-i+2}$.\\

Next we express $P$ and $Q$ in terms of the parameters
$v,~b,~k,~d$, and $i$.
As for $P$, by counting the cardinality of the following set,
$$\{(X,B)\mid B\in \mathcal B,~X\in S_1,~|X\cap B|=d+i-2\},$$
in two ways, we obtain the following equality.
\begin{eqnarray}
P=\sum_{X\in S_1}\alpha_X=b{k\choose d+i-2}{v-k\choose i-2}.
\label{equ:5.6N}
\end{eqnarray}
Note that if a 
$d+2(i-2)$ design $X$ has $d+i-2$ points in common with some 
block then $X\in S_1$. It is easy to see that 
if two distinct blocks $B$ and $C$ 
satisfy
$|X\cap B|=|X\cap C|=d+i-2$ then $|X\cap B\cap C|= d$ holds.
Therefore by counting the cardinality 
of the following set in two ways
\begin{eqnarray}
\{(X,B,C)\mid B,C\in\mathcal B,~B\neq C,~X\in S_1,~|X\cap B|=
|X\cap C|=d+i-2\},
\nonumber
\end{eqnarray}
we obtain the following equation.
\begin{eqnarray}
Q=\sum_{X\in S_1}\alpha_X(\alpha_X-1)=b\mu_d{k-d\choose i-2}^2
=b\mu_d{c\choose i-2}^2.
\nonumber
\end{eqnarray}
where $\mu_d$ is defined in Proposition \ref{pro:5.2}.
Then (\ref{equ:5.2}) implies 
\begin{eqnarray}
Q=b\frac{{v-k\choose i}{k\choose d+i-1}}{{c\choose i}
{c\choose i-1}}{c\choose i-2}^2.
\label{equ:5.7N}
\end{eqnarray}
Then (\ref{equ:5.5N}), (\ref{equ:5.6N}) and (\ref{equ:5.7N})
imply Proposition \ref{pro:5.3}.
\hfill\qed\\
\begin{pro}
Let $b=|B|$, $n=|S_1|$ and $c=k-d$.
Then we have the following formula.
\begin{eqnarray}
n=\frac{bF(v-k)}
{\prod_{l=2}^4(v-k+c-2i+l)},
\label{equ:5.8N}
\end{eqnarray}
where $F(x)$ is defined by
\begin{eqnarray}
&&F(x)={x\choose i-2}(px^2+qx+r),\nonumber\\
&&p=-(i-2){k\choose d+i-1},\\
&&q=-\frac{1}{i-1}{k\choose d+i-1}\bigg\{d^3
-(k-6i+7)d^2-\bigg((5i-6)k-(i-1)(9i-11)\bigg)d
\nonumber\\
&&\quad -(i-1)(4i-5)k+(i-1)(2i^2-2i-1)\bigg\},\\
\mbox{and}&&\nonumber\\
&&r=
\frac{1}{i-1}{k\choose d+i-1}\bigg\{(2i-3)d^3
-(3i-4)(k-3i+4)d^2
\nonumber\\
&&\quad +\bigg((i-1)k^2-(9i^2-23i+15)k+(i-1)(12i^2-33i+23)\bigg)d
\nonumber\\
&&\quad +(i-1)^2k^2-(i-1)(6i^2-16i+11)k+2(2i-3)(i-1)^3\bigg\}.
\end{eqnarray}
\label{pro:5.4}
\end{pro}
{\bf Proof}\\
For a $d+2(i-2)$ set $X$, there exists at most 
one block $B$ satisfying $|B\cap X|\geq d+i-1$.
Therefore, we have
\begin{eqnarray}&&
n={v\choose d+2i-4}-b\sum_{j=0}^{i-3}
{k\choose d+2i-4-j}{v-k\choose j}.
\label{equ:5.12}
\end{eqnarray}
Since $\mathcal D$ is a $\beta(i)$ design, the equality in (\ref{equ:1.1}) implies
the following.
\begin{eqnarray}
b=\frac{{v\choose d+2i-1}}{{k\choose d+2i-1}
+{k\choose d+2i-2}{v-k\choose 1}+\cdots+
{k\choose d+i}{v-k\choose i-1}}.
\end{eqnarray} 
Therefore, we have
\begin{eqnarray}
&&{v\choose d+2(i-2)}=
\frac{(d+2i-1)(d+2i-2)(d+2i-3)}{(v-(d+2i-4))(v-(d+2i-3))
(v-(d+2i-2))}{v\choose d+2i-1}
\nonumber\\
&&= \frac{(d+2i-1)(d+2i-2)(d+2i-3)}{(v-(d+2i-4))(v-(d+2i-3))
(v-(d+2i-2))}\times
\nonumber\\
&&b\bigg({k\choose d+2i-1}
+{k\choose d+2i-2}{v-k\choose 1}+\cdots+
{k\choose d+i}{v-k\choose i-1}\bigg).
\label{equ:5.14}
\end{eqnarray}
Then (\ref{equ:5.12}) and (\ref{equ:5.14})
imply the following equality,
\begin{eqnarray}
&&
n=\frac{(d+2i-1)(d+2i-2)(d+2i-3)}{(v-k+c-2i+4)(v-k+c-2i+3))
(v-k+c-2i+2))}\times
\nonumber\\
&&b\bigg({k\choose d+2i-1}
+{k\choose d+2i-2}{v-k\choose 1}+\cdots+
{k\choose d+i}{v-k\choose i-1}\bigg)
\nonumber\\
&&
-b\bigg({k\choose d+2i-4}+{k\choose d+2i-3}
{v-k\choose 1}+
\cdots+{k\choose d+i-1}{v-k\choose i-3}\bigg).
\nonumber\\
&&
\label{equ:5.15}
\end{eqnarray}
Therefore, if we define a polynomial $F(x)$
in $x$ by
\begin{eqnarray}
&&F(x)=\bigg\{
(d+2i-1)(d+2i-2)(d+2i-3)\times\nonumber\\
&&
\bigg(
{k\choose d+2i-1}
+{k\choose d+2i-2}{x\choose 1}+\cdots+
{k\choose d+i}{x\choose i-1}\bigg)
\nonumber\\
&&
-(x+c-2i+4)(x+c-2i+3)
(x+c-2i+2)\times\nonumber\\
&&\bigg({k\choose d+2i-4}+{k\choose d+2i-5}
{x\choose 1}+
\cdots+{k\choose d+i-1}{x\choose i-3}\bigg)
\bigg\},
\label{equ:5.16}
\end{eqnarray}
we obtain
\begin{eqnarray}
n=\frac{bF(v-k)}{(v-k+c-2i+4)(v-k+c-2i+3)
(v-k+c-2i+2)}.
\end{eqnarray}
In the following we will prove that $F(x)$ has the 
expression given in the statement of 
Proposition \ref{pro:5.4}.
For this purpose 
we define a polynomial $G(x)$ in $x$ 
of degree $i+2$ by
\begin{eqnarray}
&&G(x)= 
(d+2i-1)(d+2i-2)(d+2i-3)\times\nonumber\\
&&
 \bigg\{{k\choose d+2i-1}
 +{k\choose d+2i-2}{x\choose 1}+\cdots+{k\choose d+i}{x\choose i-1} \bigg\}
 \nonumber\\
&&
 -(x+c-2i+4)(x+c-2i+3)(x+c-2i+2)
 \times\nonumber\\
&&
\bigg\{ {k\choose d+2i-4}
 +{k\choose d+2i-5}{x\choose 1}+\cdots
 +{k\choose d+i-3}{x\choose i-1}\bigg\}.
\end{eqnarray}
Then we have
\begin{eqnarray}
&&F(x)=G(x)+(x+c-2i+4)(x+c-2i+3)(x+c-2i+2)\times\nonumber\\
&&
\left\{{c+d\choose d+i-2}{x\choose i-2}
+{c+d\choose d+i-3}{x\choose i-1}\right\}.
\label{equ:5.19}
\end{eqnarray}
We claim that the following proposition holds.
\begin{pro}
\begin{enumerate}
\item
 $G(j)=0$ for any integer $j$ with $0\leq  j \leq i-1$.
\item $F(j)=0$  for any integer $j$ with 
$0\leq j \leq i-3$.   
\item $F(i-2)=\frac{(c+d)!}{(d+i-2)!(c-i-1)!}$,

$F(i-1)=(c-i+3)(c-i+2)(c-i+1)\left\{{c+d\choose d+i-2}(i-1)
+{c+d\choose d+i-3}\right\}$.
\end{enumerate}
\label{pro:5.5}
\end{pro}

\noindent
{\bf Proof.}\\
\noindent
(1) Let $j$ be an integer satisfying $0\leq j\leq i-1$. Then by (\ref{equ:5.1}) we have
\begin{eqnarray}
&&G(j)=
(d+2i-1)(d+2i-2)(d+2i-3)
\frac{{c+d+j\choose j}{c+d\choose d+2i-1-j}}{{d+2i-1\choose j}}
\nonumber\\
&&
-(j+c-2i+4)(j+c-2i+3)(j+c-2i+2)
\frac{{c+d+j\choose j}{c+d\choose d+2i-4-j}}{{d+2i-4\choose j}}
\nonumber\\
&&={c+d+j\choose j}\bigg(
\frac{(c+d)!j!}{(d+2i-4)!
(c-2i+1+j)!}
-\frac{(c+d)!j!}
{(c-2i+1+j)!(d+2i-4)!}\bigg)\nonumber\\
&&=0.
\end{eqnarray} 
(2) For $0\leq j\leq i-3$, (1) and (\ref{equ:5.19}) imply
$F(j)=G(j)=0$.\\
(3) Since $G(i-2)=G(i-1)=0$, (\ref{equ:5.19})
implies
\begin{eqnarray}&&F(i-2)
=(c-i+2)(c-i+1)(c-i){c+d\choose d+i-2}
=\frac{(c+d)!}{(c-i-1)!(d+i-2)!}.\nonumber\\
&&F(i-1)=(c-i+3)(c-i+2)(c-i+1)
\left\{{c+d\choose d+i-2}(i-1)+{c+d\choose d+i-3}\right\}
\nonumber\\
&&=\frac{((i-1)c-i^2+d+5(i-1))(c+d)!}
{(c-i)!(d+i-2)!}.
\nonumber
\end{eqnarray}
\hfill\qed\\

\noindent
Now we are ready to finish the proof of Proposition \ref{pro:5.4}.
By (\ref{equ:5.16}), $F(x)$ is a polynomial in $x$
of degree $i$ and part (2) of Proposition \ref{pro:5.5}
shows that $0,1,\ldots,i-3$ are zeros of $F(x)$.
Hence $F(x)$ is expressed by
\begin{eqnarray}
F(x)={x\choose i-2}(px^2+qx+r),
\label{equ:5.21}
\end{eqnarray}
where $p,q$, and $r$ are some 
rational expressions of 
$c,~d$, and $i$. Moreover, since the coefficient 
of $x^i$ in $F(x)$
is $-{c+d\choose d+i-1}\frac{1}{(i-3)!}$,
we must have 
\begin{eqnarray}p=-{c+d\choose d+i-1}(i-2).
\label{equ:5.22}
\end{eqnarray}
Then, part (3) of
Proposition \ref{pro:5.5} and (\ref{equ:5.21}) imply
\begin{eqnarray}&&(c-i+2)(c-i+1)(c-i){c+d\choose d+i-2}
=p(i-2)^2+q(i-2)+r, 
\label{equ:5.23}\\
~\mbox{and}&&\nonumber\\
&&(c-i+3)(c-i+2)(c-i+1)
\left\{{c+d\choose d+i-2}(i-1)+{c+d\choose d+i-3}\right\}\nonumber\\
&&
=(i-1)(p(i-1)^2+q(i-1)+r).\label{equ:5.24}
\end{eqnarray}
Then
(\ref{equ:5.22}), (\ref{equ:5.23}) and (\ref{equ:5.24})
imply the formula of $F(x)$ given in Proposition \ref{pro:5.4}.
\hfill\qed\\

\noindent
{\bf Proof of Proposition 6}\\
Proposition \ref{pro:5.3} and (\ref{equ:5.8N}) 
imply the following inequality
\begin{eqnarray}
&&\frac{(v-k+c-2i+4)(v-k+c-2i+3)(v-k+c-2i+2)}{F(v-k)}
{c+d\choose d+i-2}{v-k\choose i-2}
\nonumber\\
&&\leq
1+\frac{(i-1)(v-k-i+1)(v-k-i+2)}{(d+i-1)(c-i+1)(c-i+2)}.
\end{eqnarray}
Then the formula of $F(x)$ given in Proposition \ref{pro:5.4} 
implies
\begin{eqnarray}
&&0\leq (p(v-k)^2+q(v-k)+r)\bigg(1
+\frac{(i-1)(v-k-i+1)(v-k-i+2)}{(d+i-1)(c-i+1)(c-i+2)}
\bigg)\nonumber
\\&&-(v-k+c-2i+4)(v-k+c-2i+3)(v-k+c-2i+2)
{c+d\choose d+i-2}.\nonumber\\
&&=\frac{(c+d)!(v-k-i+2)(cd-c+(c-d)i+d-(i-1)(v-k))}{(c-i+2)!(d+i-1)!(d+i-1)(c-i+1)(i-1)}g(v-k,d,c,i),
\nonumber\\
&&=\frac{(c+d)!(v-k-i+2)((d+2i-2)c-(i-1)v)}{(c-i+2)!(d+i-1)!(d+i-1)(c-i+1)(i-1)}g(v-k,d,c,i),
\label{equ:5.26}\end{eqnarray}
where $g(x,d,c,i)$ is defined by (\ref{equ:1.7}).
Since $v$ does not achieve the upper bound of (\ref{equ:1.6}), 
by the assumption, we have $(d+2i-2)c-(i-1)v>0$.
Hence by (\ref{equ:5.26}) we obtain $g(v-k,d,c,i)\geq 0$.

Finally, by Proposition \ref{pro:5.3} the equality 
$g(v-k,d,c,i)= 0$ holds if and only if 
 $|\{B\mid |B\cap X|=d+i-2\}|$ does not depend on the choice of $X\in S_1$.
This completes the proof of Proposition 6.
\hfill\qed\\

\noindent
{\bf Remark}
{\it
\begin{enumerate}
\item
 If $i =2$, then the inequality $g(x,d,c,i)\geq 0$ is reduced to the inequality $x\leq (cd+2)/2$.
The complementary designs of Steiner systems 
$S(t, t+1, 2t+3)$ achieve this bound.
In Theorem 2 in \cite{Noda-2001I}, 
the parameters of non trivial $\beta(2)$ 
sets are expressed by means of two
parameters and we see that there exists 
no non trivial $\beta(2)$ design 
in which $g(v-k,d,c,2)=0$ holds.
\item The parameter $(v=2k, k=d+2i-3, d)$ 
satisfies $g(v-k,d,c,i)=0$. However by Lemma {\rm\ref{lem:2.2}}  there exists no $\beta(i)$ 
design having this parameter.
\end{enumerate}
}

\noindent
{\bf Proof of Proposition 7}\\
By assumption, $v$ does not achieve the lower bound in (\ref{equ:1.6}),
i.e. $v>(d+2i)(k-d)/i$.
By Theorem 1 \cite{Noda-2001I}, the complementary design $\mathcal D'$ of $\mathcal D$ 
is a $\beta(i')$ design with the parameter 
$(v,k',d')$, where 
$$i'=k-d-i+1\geq 2,\quad k'=v-k,\quad d' = v-2k+d.$$
Since $\mathcal D$ is not a $\beta(i+1)$ design,
$\mathcal D'$ is not a
$\beta(k-d-(i+1)+1)=\beta(i'-1)$ design. 
Then Proposition 6 implies  
$g(v-k',d',c',i')=g(c+d,v-k-c,c,c-i+1)\geq 0$.
On the other hand we have
\begin{eqnarray}
&&g(c+d,x-c,c,c-i+1)
\nonumber\\
&&
=(c-i)(c-i-1)(c+d)^2
-(c-i)(c+d)\bigg(2(x-c)(2c-i)
-2c(c-i)-3c+3i+1\bigg)
\nonumber\\
&&
+(2c-i)(2c-i-1)(x-c)^2
-(2c-i)\bigg(2c(c-i+1)-6c+3i+1\bigg)(x-c)
\nonumber\\
&&
+(c-i)\bigg((c-i-1)c^2
-(3c-3i-1)c-2(2c-2i-1)(c-i)\bigg)
\nonumber\\
&&=
i(i+1)x^2-i(2(c-i)d+2ci+3i+1)x
\nonumber\\
&&
+(c-i)(c-i-1)d^2+(2ci+c-3i-1)(c-i)d+i(c^2i+c^2-3ci+4i^2-c+2i
\nonumber\\
&&=h(x,d,c,i).
\end{eqnarray}
This proves $h(v-k,d,c,i)\geq 0$.\\
Next, let $\mathcal B'$ be the block set of $\mathcal D'$. Then
$\mathcal B'=\{\Omega\backslash B\mid B\in \mathcal B\}$. 
Let 

$S_1'=\{ X' \subset \Omega \mid |X'|= d'+2(i'-2),
~|X'\cap B'|\leq d'+i'-2~ \mbox{for $\forall B' \in \mathcal B'$}\}$. \\
Then we can easily show that the following conditions are equivalent.
\begin{enumerate}
\item $X'\in S_1'$.
\item $|X'|=d'+2(i'-2)$ and $|X'\cap B'|\leq d'+i'-2 $ for any block $B'\in
\mathcal B'$.
\item $|(\Omega\backslash X)\cap (\Omega\backslash B)|
\leq v-k-i-1$ for $X=\Omega\backslash X'$ and
any block $B\in \mathcal B$.
\item $|X|=d+2(i+1)$, 
$X=\Omega\backslash X'$ and $|X\cap B|\leq d+i+1$ for
any $B\in\mathcal B$.
\item $X\in S_2$.
\end{enumerate}
Also the following conditions are equivalent.
\begin{enumerate}
\item[(6)] $|X'\cap B'|= d'+i'-2$ for $X'\in S_1'$ and a block
$B'\in \mathcal B'$. 
\item[(7)] $|X\cap B|= d+i+1$ for $X\in S_2$
and a block $B\in \mathcal B$.
\end{enumerate}
Therefore Proposition 6 implies $h(v-k,d,c,i)=0$
if and only if $|\{B\in \mathcal B\mid |X\cap B|=
d+i-1\}|$ does not depend on the choice of $X\in S_2$.
This completes the proof of Proposition 7.
\hfill\qed\\

\noindent
Finally we present a proof of Lemma \ref{lem:5.1}.\\

\noindent
{\bf Proof of Lemma \ref{lem:5.1}}.\\
\noindent
First we assume $m\geq s$. Let $\Omega_1$ be a
$(m+j)$ point set. We count the cardinality of the following
set in two ways.
\begin{eqnarray}\{(X,Y)\mid X,~Y\subset \Omega_1,
|X|=s,~|Y|=j,~X\cap Y=\emptyset\}.
\label{equ:5.30}
\end{eqnarray}
For each $Y\subset \Omega_1$ with $|Y|=j$
there is ${m\choose s}$ choices for $X$. 
Therefore the cardinality of the set given in
(\ref{equ:5.30})
equals ${m+j\choose j}{m\choose s}$.
Next, we take a $j$ point subset $W$ of $\Omega_1$
and fix it. For a $s$ point set $X\subset \Omega_1$, we have
$0\leq |X\cap W|\leq j$. The number of $s$ point subsets
$X$ with $|X\cap W|=l$ is ${j\choose l}{m\choose s-l}$
for $0\leq l\leq j$ and each $X$ with $|X\cap W|=l$
the number of the choices of $Y$ is ${m+j-s\choose j}$.
Hence the cardinality of the set in (\ref{equ:5.30})
equals
\begin{eqnarray}
\bigg({m\choose s}+{m\choose s-1}{j\choose 1}
+{m\choose s-2}{j\choose 2}+\cdots
+{m\choose s-j}{j\choose j}\bigg){m+j-s\choose j}.
\end{eqnarray}
Since ${m+j-s\choose j}>0$, we have
\begin{eqnarray}
&&{m\choose s}+{m\choose s-1}{j\choose 1}
+{m\choose s-2}{j\choose 2}+\cdots
+{m\choose s-j}{j\choose j}
=\frac{{m+j\choose j}{m\choose s}}{{m+j-s\choose j}
}
\nonumber\\
&&=\frac{{m+j\choose j}{m\choose s-j}}{{s\choose j}}.
\end{eqnarray}
Next we assume $s>m\geq s-j$. Since
${m\choose s}={m\choose s-1}=\cdots={m\choose m+1}=0$,
the left side of (\ref{equ:5.1}) equals
\begin{eqnarray}&&{m\choose m}{j\choose s-m}+{m\choose m-1}{j\choose s-m+1}
+\cdots+{m\choose s-j}{j\choose j}
\nonumber\\
&&={m\choose 0}{j\choose m+j-s}
+{m\choose 1}{j\choose (m+j-s)-1}
+\cdots+{m\choose m+j-s}{j\choose 0}.
\nonumber\\
&&
\label{equ:5.33}
\end{eqnarray}
By the assumption we have $0\leq m+j-s\leq j$.
Then we consider the 
following set consisting of pairs of subset in $(m+j)$
point set $\Omega_1$.
\begin{eqnarray}&&
\{(X,Y)\mid  X,Y\subset \Omega_1,~|X|=m+j-s,~|Y|=j,~
X\cap Y=\emptyset\}.
\label{equ:5.34}
\end{eqnarray}
We count the cardinality of the set defined 
in (\ref{equ:5.34}) in two ways. 
Then by a similar consideration we have
\begin{eqnarray}
{m\choose 0}{j\choose m+j-s}+{m\choose 1}{j\choose m+j-s-1}
+\cdots+{m\choose m+j-s}{j\choose 0}
=\frac{{m+j\choose j}{m\choose m+j-s}}{{s\choose j}}.
\nonumber
\end{eqnarray}
By (\ref{equ:5.33}) this implies (\ref{equ:5.1}).
Finally if $s>m+j$, then both side of
(\ref{equ:5.1}) equals 0.
This completes the proof of Lemma \ref{lem:5.1}.
\hfill \qed

\section{Some Open Problems}
Finally we present some open problems in $\beta(i)$ designs. 
The most fundamental problem is the following.\\

\noindent
{\bf Problem 1}. {\it Does there exist a $\beta(i)$ design which has the parameter 
$(v,k,d)$ with $d>0$ and
$1<i<k-d$ other than Steiner system $S(5,8,24)$ and the complementary design of 
it {\rm ?}}\\

The Steiner system $S(5,8,24)$ is both a $\beta(1)$ and a $\beta(2)$ design and the complementary design of it is both a $\beta(3)$ and a $\beta(4)$ design. 
Hauck \cite{Hauck-1982} proved that if $\mathcal D$ is both a $\beta(1)$ and 
a $\beta(2)$ design then $\mathcal D$ is the Steiner system $S(5,8,24)$ or a 
Steiner system $S(t,t+1,2t+2)$ for some $t$. 
So our second problem is the following.\\

\noindent
{\bf Problem 2.} {\it Does there exist a $\beta(i)$ design which is also a 
$\beta(i+1)$ design other than the Steiner system $S(5,8,24)$,
the complementary design of it and Steiner systems $S(t,t+1,2t+2)$ {\rm ?}}\\

 It is proved in \cite{Noda-2001I} that there exists no $\beta(i)$ design which is also a 
$\beta(i+1)$ design if $i\equiv 2~(\mbox{mod}~ 4)$. 
Therefore there exists no $\beta(2)$ design which is also a $\beta(3)$ design. In a forthcoming paper \cite{Bannai-Noda} we prove that if $\mathcal D$ is both a 
$\beta(3)$ and a $\beta(4)$ design with the parameter $(v,k,d)$,
then $\mathcal D$ is the complementary design of the Steiner system 
$S(5,8,24)$ or possibly $v=2k$ and $k=d+6$. Generally $\beta(i)$ designs with the parameter 
$(v=2k, k=d+2i,d)$ achieve the lower bound of (\ref{equ:1.6}),
whence are also $\beta(i+1)$ designs. In the case $i=1$ 
such designs do exist. Our third problem is the following.\\

\noindent
{\bf Problem 3}. {\it Does there exist a $\beta(i)$ design which has the 
parameter $(v=2k, k=d+2i,d)$ for $i>1$ {\rm ?}}\\

 Whether there exists a perfect $e$-code $\mathcal C$ in $J(v,k)$ 
with $|\mathcal C|\geq 2$  and $1<d(\mathcal C)<k$ is 
a well known open problem. 
This is stated in terms of a $\beta(i)$ design as follows.\\

\noindent
{\bf Problem 4}. {\it Does there exist a $\beta(i)$ design which has the parameter 
$(v,k=d+2i-1,d)$ with $0< d< k-1$ {\rm ?}}\\

\noindent
{\bf Acknowledgments} \\
 The authors thank Eiichi Bannai for valuable suggestions
 especially telling us the
 relation between $\beta(i)$ designs and 
 perfect codes in $J(v,k)$, and introducing the work of Etzion-Schwartz
 \cite{Etzion-S-2004} and 
Ahlswede-Aydinian-Khachatrian \cite{Ahlswede-A-K-2001}.
We also would like to express our gratitude to the referee for helpful comments and
advice.

\noindent
Etsuko Bannai: Misakigaoka 2-8-21, Itoshima-shi, Fukuoka, 819-1136, Japan\\
e-mail: et-ban@rc4.so-net.ne.jp 
 \\

\noindent
Ryuzaburo Noda: Yokoikami 507-79, Kita-ku, Okayama, 701-1145 , Japan\\
e-mail: rrnoda@yahoo.co.jp

\end{document}